\documentclass[12pt,twoside]{article}
\usepackage[english]{babel}
\usepackage[latin1]{inputenc}
\usepackage{amsmath}
\usepackage{graphicx}
\usepackage{times,amssymb,amscd}
\newtheorem{thm}{Theorem}[section]

\newtheorem{lem}[thm]{Lemma}

\newtheorem{defn}[thm]{Definition}
\newtheorem{rem}[thm]{Remark}
\numberwithin{equation}{section}

\newcommand{\bA}{\mathbf{A}}

\newcommand{\bE}{\mathbf{E}}

\newcommand{\bH}{\mathbf{H}}

\newcommand{\bL}{\mathbf{L}}

\newcommand{\bR}{\mathbf{R}}
\newcommand{\bS}{\mathbf{S}}
\newcommand{\bV}{\mathbf{V}}

\newcommand{\be}{\mathbf{e}}

\newcommand{\bT}{\mathbf{T}}

\newcommand{\bu}{\mathbf{u}}
\newcommand{\bv}{\mathbf{v}}
\newcommand{\bt}{\mathbf{t}}

\newcommand{\EUC}{\mathbf E^3}
\newcommand{\SPH}{\bS^3}
\newcommand{\HYP}{\bH^3}
\newcommand{\SXR}{\bS^2\!\times\!\bR}
\newcommand{\HXR}{\bH^2\!\times\!\bR}
\newcommand{\SLR}{\widetilde{\bS\bL_2\bR}}
\newcommand{\NIL}{\mathbf{Nil}}
\newcommand{\SOL}{\mathbf{Sol}}

\begin{document}
\pagestyle{myheadings}
\markboth{\centerline{Jen\H o Szirmai}}
{Triangle angle sums related to translation curves $\dots$}
\title
{Triangle angle sums related to translation curves in $\SOL$ geometry
\footnote{Mathematics Subject Classification 2010: 53A20, 53A35, 52C35, 53B20. \newline
Key words and phrases: Thurston geometries, $\SOL$ geometry, translation and geodesic triangles, interior angle sum \newline
}}

\author{Jen\H o Szirmai \\
\normalsize Budapest University of Technology and \\
\normalsize Economics Institute of Mathematics, \\
\normalsize Department of Geometry \\
\normalsize Budapest, P. O. Box: 91, H-1521 \\
\normalsize szirmai@math.bme.hu
\date{\normalsize{\today}}}

\maketitle
\begin{abstract}
After having investigated the geodesic and translation triangles and their angle sums in $\NIL$ and $\SLR$ geometries we consider the analogous problem in $\SOL$ space that
is one of the eight 3-dimensional Thurston geometries.

We analyse the interior angle sums of translation triangles in $\SOL$ geometry
and prove that it can be larger or equal than $\pi$.

In our work we will use the projective model of $\SOL$ described by E. Moln\'ar in \cite{M97},
\end{abstract}

\section{Introduction} \label{section1}
In the Thurston spaces can be introduced in a natural way (see \cite{M97}) translations mapping each point to any point.
Consider a unit vector at the origin. Translations, postulated at the
beginning carry this vector to any point by its tangent mapping. If a curve $t\rightarrow (x(t),y(t),z(t))$ has just the translated
vector as tangent vector in each point, then the  curve is
called a {\it translation curve}. This assumption leads to a system of first order differential equations, thus translation
curves are simpler than geodesics and differ from them in $\NIL$, $\SLR$ and $\SOL$ geometries. In $\EUC$, $\SPH$, $\HYP$, $\SXR$ and $\HXR$ geometries the mentioned curves
coincide with each other.

Therefore, the translation curves also play an important role in $\NIL$, $\SLR$ and $\SOL$ geometries and often seem to be more natural in these geometries,
than their geodesic lines.

A translation triangle in Riemannian geometry and more generally in metric geometry a
figure consisting of three different points together with the pairwise-connecting translation curves.
The points are known as the vertices, while the translation curve segments are known as the sides of the triangle.

In the geometries of constant curvature $\EUC$, $\HYP$, $\SPH$ the well-known sums of the interior angles of geodesic (or translation)
triangles characterize the space. It is related to the Gauss-Bonnet theorem which states that the integral of the Gauss curvature
on a compact $2$-dimensional Riemannian manifold $M$ is equal to $2\pi\chi(M)$ where $\chi(M)$ denotes the Euler characteristic of $M$.
This theorem has a generalization to any compact even-dimensional Riemannian manifold (see e.g. \cite{Ch}, \cite{KN}).

In \cite{CsSz16} we investigated the angle sum of translation and geodesic triangles in $\SLR$ geometry
and proved that the possible sum of the interior angles in
a translation triangle must be greater or equal than $\pi$. However, in geodesic triangles this sum is
less, greater or equal to $\pi$.

In \cite{Sz16} we considered the analogous problem for geodesic triangles in $\NIL$ geometry and proved
that the sum of the interior angles of geodesic triangles in $\NIL$ space is larger. less or equal than $\pi$.
In \cite{B} K.~Brodaczewska showed, that sum of the interior angles of translation triangles of the $\NIL$ space is larger than $\pi$.

However, in $\SXR$, $\HXR$ and $\SOL$ Thur\-ston geo\-metries there are no result concerning the
angle sums of translation or geodesic triangles. Therefore, it is interesting to study similar question
in the above three geometries.
Now, we are interested in {\it translation triangles} in $\SOL$ space \cite{S,T}.

In Section 2 we describe the projective model and the isometry group of $\SOL$,
moreover, we give an overview about its translation curves.
\begin{rem}
We note here, that nowadays the $\SOL$ geometry is a widely investigated space concerning
its manifolds, tilings, geodesic and translation ball packings and probability theory
(see e.g. \cite{BT}, \cite{CaMoSpSz}, \cite{KV}, \cite{MSz}, \cite{MSz12}, \cite{MSzV}, \cite{Sz13-2} and the references given there).
\end{rem}
{\it In Section 3 we study the $\SOL$ translation triangles and prove that their interior angle sums can be larger or equal than $\pi$.}
\section{On Sol geometry}
\label{sec:1}

In this Section we summarize the significant notions and notations of real $\SOL$ geometry (see \cite{M97}, \cite{S}).

$\SOL$ is defined as a 3-dimensional real Lie group with multiplication
\begin{equation}
     \begin{gathered}
(a,b,c)(x,y,z)=(x + a e^{-z},y + b e^z ,z + c).
     \end{gathered} \tag{2.1}
     \end{equation}
We note that the conjugacy by $(x,y,z)$ leaves invariant the plane $(a,b,c)$ with fixed $c$:
\begin{equation}
     \begin{gathered}
(x,y,z)^{-1}(a,b,c)(x,y,z)=(x(1-e^{-c})+a e^{-z},y(1-e^c)+b e^z ,c).
     \end{gathered} \tag{2.2}
     \end{equation}
Moreover, for $c=0$, the action of $(x,y,z)$ is only by its $z$-component, where $(x,y,z)^{-1}=(-x e^{z}, -y e^{-z} ,-z)$. Thus the $(a,b,0)$ plane is distinguished as a {\it base plane} in
$\SOL$, or by other words, $(x,y,0)$ is normal subgroup of $\SOL$.
$\SOL$ multiplication can also be affinely (projectively) interpreted by "right translations"
on its points as the following matrix formula shows, according to (2.1):
     \begin{equation}
     \begin{gathered}
     (1;a,b,c) \to (1;a,b,c)
     \begin{pmatrix}
         1&x&y&z \\
         0&e^{-z}&0&0 \\
         0&0&e^z&0 \\
         0&0&0&1 \\
       \end{pmatrix}
       =(1;x + a e^{-z},y + b e^z ,z + c)
       \end{gathered} \tag{2.3}
     \end{equation}
by row-column multiplication.
This defines "translations" $\mathbf{L}(\mathbf{R})= \{(x,y,z): x,~y,~z\in \mathbf{R} \}$
on the points of space $\SOL= \{(a,b,c):a,~b,~c \in \mathbf{R}\}$.
These translations are not commutative, in general.
Here we can consider $\mathbf{L}$ as projective collineation group with right actions in homogeneous
coordinates as usual in classical affine-projective geometry.
We will use the Cartesian homogeneous coordinate simplex $E_0(\be_0)$, $E_1^{\infty}(\be_1)$, \ $E_2^{\infty}(\be_2)$, \ 
$E_3^{\infty}(\be_3), \ (\{\be_i\}\subset \bV^4$ \ $\text{with the unit point}$ $E(\be = \be_0 + \be_1 + \be_2 + \be_3 ))$
which is distinguished by an origin $E_0$ and by the ideal points of coordinate axes, respectively.
Thus {$\SOL$} can be visualized in the affine 3-space $\bA^3$
(so in Euclidean space $\bE^3$) as well \cite{M97}.

In this affine-projective context E. Moln\'ar has derived in \cite{M97} the usual infinitesimal arc-length square at any point
of $\SOL$, by pull back translation, as follows
\begin{equation}
   \begin{gathered}
      (ds)^2:=e^{2z}(dx)^2 +e^{-2z}(dy)^2 +(dz)^2.
       \end{gathered} \tag{2.4}
     \end{equation}
Hence we get infinitesimal Riemann metric invariant under translations, by the symmetric metric tensor field $g$ on $\SOL$ by components as usual.

It will be important for us that the full isometry group Isom$(\SOL)$ has eight components, since the stabilizer of the origin
is isomorphic to the dihedral group $\mathbf{D_4}$, generated by two involutive (involutory) transformations, preserving (2.4):
\begin{equation}
   \begin{gathered}
      (1)  \ \ y \leftrightarrow -y; \ \ (2)  \ x \leftrightarrow y; \ \ z \leftrightarrow -z; \ \ \text{i.e. first by $3\times 3$ matrices}:\\
     (1) \ \begin{pmatrix}
               1&0&0 \\
               0&-1&0 \\
               0&0&1 \\
     \end{pmatrix}; \ \ \
     (2) \ \begin{pmatrix}
               0&1&0 \\
               1&0&0 \\
               0&0&-1 \\
     \end{pmatrix}; \\
     \end{gathered} \tag{2.5}
     \end{equation}
     with its product, generating a cyclic group $\mathbf{C_4}$ of order 4
     \begin{equation}
     \begin{gathered}
     \begin{pmatrix}
                    0&1&0 \\
                    -1&0&0 \\
                    0&0&-1 \\
     \end{pmatrix};\ \
     \begin{pmatrix}
               -1&0&0 \\
               0&-1&0 \\
               0&0&1 \\
     \end{pmatrix}; \ \
     \begin{pmatrix}
               0&-1&0 \\
               1&0&0 \\
               0&0&-1 \\
     \end{pmatrix};\ \
     \mathbf{Id}=\begin{pmatrix}
               1&0&0 \\
               0&1&0 \\
               0&0&1 \\
     \end{pmatrix}.
     \end{gathered} \notag
     \end{equation}
     Or we write by collineations fixing the origin $O(1,0,0,0)$:
\begin{equation}
(1) \ \begin{pmatrix}
         1&0&0&0 \\
         0&1&0&0 \\
         0&0&-1&0 \\
         0&0&0&1 \\
       \end{pmatrix}, \ \
(2) \ \begin{pmatrix}
         1&0&0&0 \\
         0&0&1&0 \\
         0&1&0&0 \\
         0&0&0&-1 \\
       \end{pmatrix} \ \ \text{of form (2.3)}. \tag{2.6}
\end{equation}
A general isometry of $\SOL$ to the origin $O$ is defined by a product $\gamma_O \tau_X$, first $\gamma_O$ of form (2.6) then $\tau_X$ of (2.3). To
a general point $A(1,a,b,c)$, this will be a product $\tau_A^{-1} \gamma_O \tau_X$, mapping $A$ into $X(1,x,y,z)$.

Conjugacy of translation $\tau$ by an above isometry $\gamma$, as $\tau^{\gamma}=\gamma^{-1}\tau\gamma$ also denotes it, will also be used by
(2.3) and (2.6) or also by coordinates with above conventions.

We remark only that the role of $x$ and $y$ can be exchanged throughout the paper, but this leads to the mirror interpretation of $\SOL$.
As formula (2.4) fixes the metric of $\SOL$, the change above is not an isometry of a fixed $\SOL$ interpretation. Other conventions are also accepted
and used in the literature.

{\it $\SOL$ is an affine metric space (affine-projective one in the sense of the unified formulation of \cite{M97}). Therefore its linear, affine, unimodular,
etc. transformations are defined as those of the embedding affine space.}
\subsection{Translation curves}

We consider a $\SOL$ curve $(1,x(t), y(t), z(t) )$ with a given starting tangent vector at the origin $O(1,0,0,0)$
\begin{equation}
   \begin{gathered}
      u=\dot{x}(0),\ v=\dot{y}(0), \ w=\dot{z}(0).
       \end{gathered} \tag{2.7}
     \end{equation}
For a translation curve let its tangent vector at the point $(1,x(t), y(t), z(t) )$ be defined by the matrix (2.3)
with the following equation:
\begin{equation}
     \begin{gathered}
     (0,u,v,w)
     \begin{pmatrix}
         1&x(t)&y(t)&z(t) \\
         0&e^{-z(t)}&0&0 \\
         0&0&e^{z(t)}& 0 \\
         0&0&0&1 \\
       \end{pmatrix}
       =(0,\dot{x}(t),\dot{y}(t),\dot{z}(t)).
       \end{gathered} \tag{2.8}
     \end{equation}
Thus, {\it translation curves} in $\SOL$ geometry (see \cite{MoSzi10} and \cite{MSz}) are defined by the first order differential equation system
$\dot{x}(t)=u e^{-z(t)}, \ \dot{y}(t)=v e^{z(t)},  \ \dot{z}(t)=w,$ whose solution is the following:
\begin{equation}
   \begin{gathered}
     x(t)=-\frac{u}{w} (e^{-wt}-1), \ y(t)=\frac{v}{w} (e^{wt}-1),  \ z(t)=wt, \ \mathrm{if} \ w \ne 0 \ \mathrm{and} \\
     x(t)=u t, \ y(t)=v t,  \ z(t)=z(0)=0 \ \ \mathrm{if} \ w =0.
       \end{gathered} \tag{2.9}
\end{equation}

We assume that the starting point of a translation curve is the origin, because we can transform a curve into an
arbitrary starting point by translation (2.3), moreover, unit velocity translation can be assumed :
\begin{equation}
\begin{gathered}
        x(0)=y(0)=z(0)=0; \\ \ u=\dot{x}(0)=\cos{\theta} \cos{\phi}, \ \ v=\dot{y}(0)=\cos{\theta} \sin{\phi}, \ \ w=\dot{z}(0)=\sin{\theta}; \\
        - \pi \leq \phi \leq \pi, \ -\frac{\pi}{2} \leq \theta \leq \frac{\pi}{2}. \tag{2.10}
\end{gathered}
\end{equation}
\begin{defn}
The translation distance $d^t(P_1,P_2)$ between the points $P_1$ and $P_2$ is defined by the arc length of the above translation curve
from $P_1$ to $P_2$.
\end{defn}
Thus we obtain the parametric equation of the the {\it translation curve segment} $t(\phi,\theta,t)$ with starting point at the origin in direction
\begin{equation}
\bt(\phi, \theta)=(\cos{\theta} \cos{\phi}, \cos{\theta} \sin{\phi}, \sin{\theta}) \tag{2.11}
\end{equation}
where $t \in [0,r \in \bR^+$]. If $\theta \ne 0$ then the system of equation is:
\begin{equation}
\begin{gathered}
        \left\{ \begin{array}{ll}
        x(\phi,\theta,t)=-\cot{\theta} \cos{\phi} (e^{-t \sin{\theta}}-1), \\
        y(\phi,\theta,t)=\cot{\theta} \sin{\phi} (e^{t \sin{\theta}}-1), \\
        z(\phi,\theta,t)=t \sin{\theta}.
        \end{array} \right. \\
        \text{If $\theta=0$ then}: ~  x(t)=t\cos{\phi} , \ y(t)=t \sin{\phi},  \ z(t)=0.
        \tag{2.12}
\end{gathered}
\end{equation}
\section{Translation triangles}
We consider $3$ points $A_1$, $A_2$, $A_3$ in the projective model of $\SOL$ space (see Section 2).
The {\it translation segments} $a_k$ connecting the points $A_i$ and $A_j$
$(i<j,~i,j,k \in \{1,2,3\}, k \ne i,j$) are called sides of the {\it translation triangle} with vertices $A_1$, $A_2$, $A_3$.
\begin{figure}[ht]
\centering
\includegraphics[width=14cm]{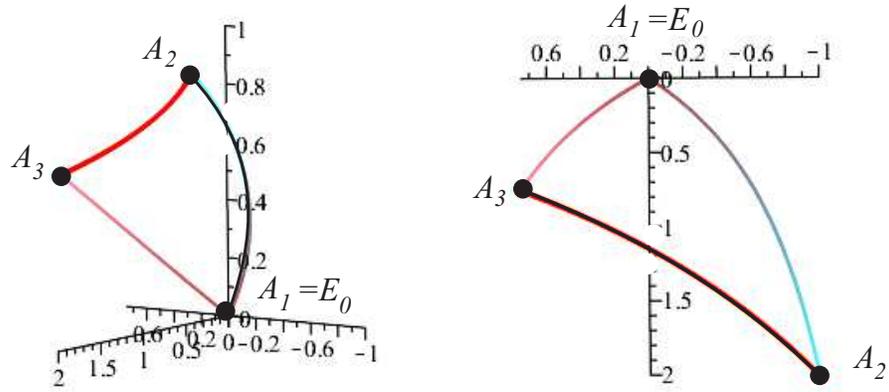}
\caption{Translation triangle with vertices $A_1=(1,0,0,0)$, $A_2=(1,-1,2,1)$, $A_3=(1,3/4,3/4,1/2)$.}
\label{}
\end{figure}
In Riemannian geometries the metric tensor (or infinitesimal arc-lenght square (see (2.4)) is used to define the angle $\theta$ between two geodesic curves.
If their tangent vectors in their common point are $\bu$ and $\bv$ and $g_{ij}$ are the components of the metric tensor then
\begin{equation}
\cos(\theta)=\frac{u^i g_{ij} v^j}{\sqrt{u^i g_{ij} u^j~ v^i g_{ij} v^j}} \tag{3.1}
\end{equation}
It is clear by the above definition of the angles and by the infinitesimal arc-lenght square (2.4), that
the angles are the same as the Euclidean ones at the origin of
the projective model of $\SOL$ geometry.

Considering a translation triangle $A_1A_2A_3$ we can assume by the homogeneity of the $\SOL$ geometry that one of its vertex
coincide with the origin $A_1=E_0=(1,0,0,0)$ and the other two vertices are $A_2(1,x^2,y^2,z^2)$ and $A_3(1,x^3,y^3,z^3)$.

We will consider the {\it interior angles} of translation triangles that are denoted at the vertex $A_i$ by $\omega_i$ $(i\in\{1,2,3\})$.
We note here that the angle of two intersecting translation curves depends on the orientation of their tangent vectors. 

{\it In order to determine the interior angles of a translation triangle $A_1A_2A_3$
and its interior angle sum $\sum_{i=1}^3(\omega_i)$,
we define {translations} $\bT_{A_i}$, $(i\in \{2,3\})$ as elements of the isometry group of $\SOL$, that
maps the origin $E_0$ onto $A_i$} (see Fig.~2).

E.g. the isometrie $\bT_{A_2}$ and its inverse (up to a positive determinant factor) can be given by:
\begin{equation}
\bT_{A_2}=
\begin{pmatrix}
1 & x^2 & y^2 & z^2 \\
0 & \mathrm{e}^{-z^2} & 0 & 0 \\
0 & 0 & \mathrm{e}^{z^2} & 0 \\
0 & 0 & 0 & 1
\end{pmatrix} , ~ ~ ~
\bT_{A_2}^{-1}=
\begin{pmatrix}
1 & -x^2\mathrm{e}^{z^2} & -y^2\mathrm{e}^{-z^2} & -z^2 \\
0 & \mathrm{e}^{z^2} & 0 & 0 \\
0 & 0 & \mathrm{e}^{-z^2} & 0 \\
0 & 0 & 0 & 1
\end{pmatrix} , \tag{3.2}
\end{equation}
and the images $\bT^{-1}_{A_2}(A_i)$ of the vertices $A_i$ $(i \in \{1,2,3\})$ are the following (see also Fig.~2):
\begin{equation}
\begin{gathered}
\bT^{-1}_{A_2}(A_1)=A_1^2=(1,-x^2\mathrm{e}^{z^2},-y^2\mathrm{e}^{-z^2},-z^2),~\bT^{-1}_{A_2}(A_2)=A_2^2=E_0=(1,0,0,0), \\ \bT^{-1}_{A_2}(A_3)=
A_3^2=(1,(x^3-x^2)\mathrm{e}^{z^2},(y^3-y^2)\mathrm{e}^{-z^2},z^3-z^2). \tag{3.3}
\end{gathered}
\end{equation}
Similarly to the above computation we get that the images $\bT^{-1}_{A_3}(A_i)$ of the vertices $A_i$ $(i \in \{1,2,3\})$ are the following (see also Fig.~2):
\begin{equation}
\begin{gathered}
\bT^{-1}_{A_3}(A_1)=A_1^3=(1,-x^3\mathrm{e}^{z^3},-y^3\mathrm{e}^{-z^3},-z^3),~\bT^{-1}_{A_3}(A_3)=A_2^2=E_0=(1,0,0,0), \\ \bT^{-1}_{A_3}(A_2)=
A_2^3=(1,(x^2-x^3)\mathrm{e}^{z^3},(y^2-y^3)\mathrm{e}^{-z^3},z^2-z^3). \tag{3.4}
\end{gathered}
\end{equation}
\begin{figure}[ht]
\centering
\includegraphics[width=13cm]{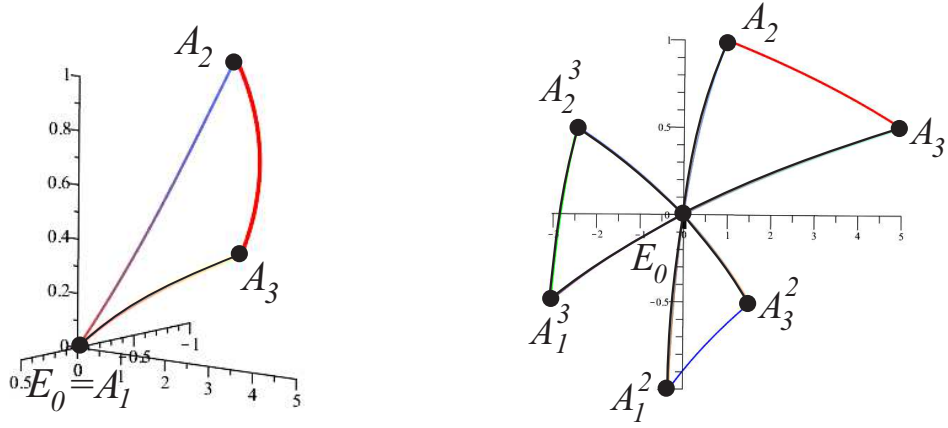}
\caption{Translation triangle with vertices $A_1=(1,0,0,0)$, $A_2=(1,-1,1,1)$, $A_3=(1,1/2,5,1/2)$ and its translated copies $A_1^2A_3^2E_0$ and $A_1^3A_2^3E_0$.}
\label{}
\end{figure}
Our aim is to determine angle sum $\sum_{i=1}^3(\omega_i)$ of the interior angles of translation triangles $A_1A_2A_3$ (see Fig.~1-2).
We have seen that $\omega_1$ and the angle of translation curves with common point at the origin $E_0$ is the same as the
Euclidean one therefore can be determined by usual Euclidean sense.

The translations $\bT_{A_i}$ $(i=2,3)$ are isometries in $\SOL$ geometry thus
$\omega_i$ is equal to the angle $(t(A_i^i, A_1^i)t(A_i^i, A_j^i))\angle$ $(i,j=2,3$, $i \ne j)$ (see Fig.~2)
where $t(A_i^i, A_1^i)$, $t(A_i^i, A_j^i)$ are oriented translation curves $(E_0=A_2^2=A_3^3)$ and
$\omega_1$ is equal to the angle $(t(E_0, A_2)t(E_0, A_3)) \angle$
where $t(E_0, A_2)$, $t(E_0, A_3)$ are also oriented translation curves.

We denote the oriented unit tangent vectors of the oriented geodesic curves $t(E_0, A_i^j)$ with $\mathbf{t}_i^j$ where
$(i,j)\in\{(1,3),(1,2),(2,3),(3,2),(3,0),(2,0)\}$ and $A_3^0=A_3$, $A_2^0=A_2$.
The Euclidean coordinates of $\mathbf{t}_i^j$ (see Section 2.1) are :
\begin{equation}
\mathbf{t}_i^j=(\cos(\theta_i^j) \cos(\alpha_i^j), \cos(\theta_i^j) \sin(\alpha_i^j), \sin(\theta_i^j)). \tag{3.5}
\end{equation}
In order to obtain the angle of two translation curves $t_{E_0A_i^j}$ and $t_{E_0A_k^l}$ ($(i,j)\ne(k,l)$; $(i,j),(k,l)\in\{(1,3),(1,2),(2,3),(3,2),(3,0),(2,0)\})$ intersected at the origin $E_0$ we need to determine their tangent vectors $\bt_s^r$
$((s,r) \in \{(1,3),(1,2),$ $(2,3),(3,2),(3,0),(2,0)\})$
(see (3.5)) at their starting point $E_0$.
From (3.5) follows that a tangent vector at the origin is given by the parameters $\phi$ and $\theta$ of the corresponding translation curve (see (2.12)) that
can be determined from the homogeneous coordinates of the endpoint of the translation curve as the following Lemma shows:
\begin{lem}
\begin{enumerate}
\item Let $(1,x,y,z)$ $(y,z \in \bR \setminus \{0\}, x \in \bR)$ be the homogeneous coordinates of the point $P \in \SOL$. The paramerters of the
corresponding translation curve $t_{E_0P}$ are the following
\begin{equation}
\begin{gathered}
\phi=\mathrm{arccot}\Big(-\frac{x}{y} \frac{\mathrm{e}^z-1}{\mathrm{e}^{-z}-1}\Big),~\theta=\mathrm{arccot}\Big( \frac{y}{\sin\phi(\mathrm{e}^z-1)}\Big),\\
t=\frac{z}{\sin\theta}, ~ \text{where} ~ -\pi < \phi \le \pi, ~ -\pi/2\le \theta \le \pi/2, ~ t\in \bR^+.
\end{gathered} \tag{3.6}
\end{equation}
\item Let $(1,x,0,z)$ $(x,z \in \bR \setminus \{0\})$ be the homogeneous coordinates of the point $P \in \SOL$. The paramerters of the
corresponding translation curve $t_{E_0P}$ are the following
\begin{equation}
\begin{gathered}
\phi=0~\text{or}~  \pi, ~\theta=\mathrm{arccot}\Big( \mp \frac{x}{(\mathrm{e}^{-z}-1)}\Big),\\
t=\frac{z}{\sin\theta}, ~ \text{where}  ~ -\pi/2\le \theta \le \pi/2, ~ t\in \bR^+.
\end{gathered} \tag{3.7}
\end{equation}
\item Let $(1,x,y,0)$ $(x,y \in \bR)$ be the homogeneous coordinates of the point $P \in \SOL$. The paramerters of the 
corresponding translation curve $t_{E_0P}$ are the following
\begin{equation}
\begin{gathered}
\phi=\arccos\Big(\frac{x}{\sqrt{x^2+y^2}}\Big),~  \theta=0,\\
t=\sqrt{x^2+y^2}, ~ \text{where}  ~ -\pi < \phi \le \pi, ~ t\in \bR^+.
\end{gathered} \tag{3.8}
\end{equation}
\end{enumerate}
\end{lem}
\begin{thm}
The sum of the interior angles of a translation triangle is greather or equal to $\pi$.
\end{thm}
\textbf{Proof:} The translations $\bT_{A_2}^{-1}$ and $\bT_{A_3}^{-1}$ are isometries
in $\SOL$ geometry thus $\omega_2$ is equal to the angle $((A_2^2 A_1^2), (A_2^2 A_3^2)) \angle$ (see Fig.~2)
of the oriented translation segments $t_{A_2^2 A_1^2}$, $t_{A_2^2A_3^2}$ and $\omega_3$ is equal to the angle
$((A_3^3 A_1^3),(A_3^3 A_2^3)) \angle$
of the oriented translation segments $t_{A_3^3 A_1^3}$ and $t_{A_3^3 A_2^3}$ $(E_0=A_2^2=A_3^3$).

Substituting the coordinates of the points $A_i^j$ (see (3.3) and (3.4)) $((i,j) \in \{(1,3),(1,2),$ $(2,3),(3,2),(3,0),(2,0)\})$ to the appropriate equations of Lemma 3.1,
it is easy to see that
\begin{equation}
\begin{gathered}
\theta_2^0=-\theta_1^2,~\phi_2^0-\phi_1^2=\pm \pi \Rightarrow \bt_2^0=-\bt_1^2,\\
\theta_3^0=-\theta_1^3,~\phi_3^0-\phi_1^3=\pm \pi \Rightarrow \bt_3^0=-\bt_1^3,\\
\theta_3^2=-\theta_2^3,~\phi_3^2-\phi_2^3=\pm \pi \Rightarrow \bt_3^2=-\bt_2^3.
\end{gathered} \tag{3.9}
\end{equation}
\begin{figure}[ht]
\centering
\includegraphics[width=8cm]{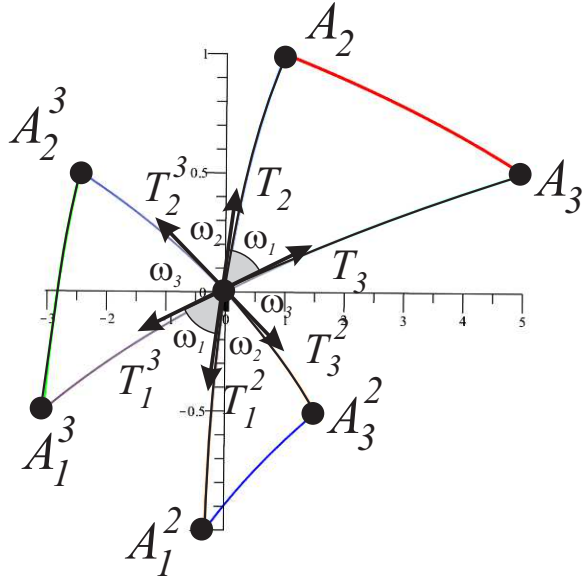}
\caption{Translation triangle with vertices $A_1=(1,0,0,0)$, $A_2=(1,-1,1,1)$, $A_3=(1,1/2,5,1/2)$ and its translated copies $A_1^2A_3^2E_0$ and $A_1^3A_2^3E_0$.}
\label{}
\end{figure}
\begin{figure}
\centering
\includegraphics[width=8cm]{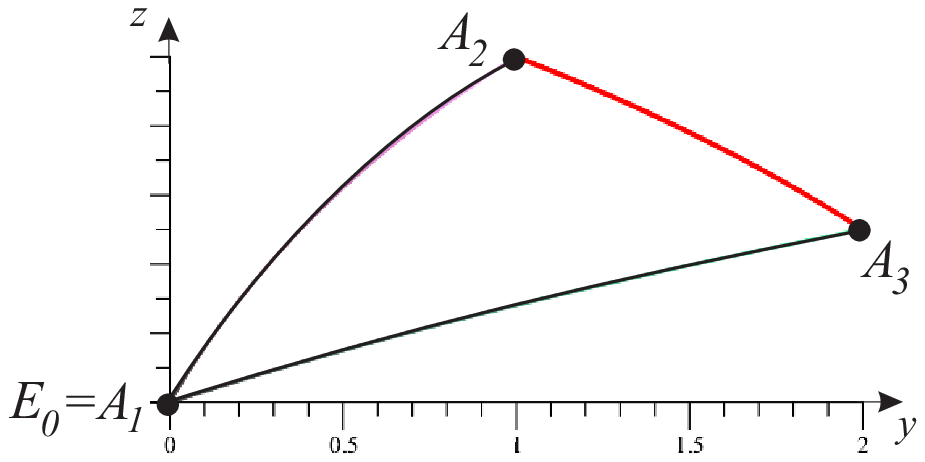}
\caption{Translation triangle with vertices $A_1=(1,0,0,0)$, $A_2=(1,0,1,1)$, $A_3=(1,0,2,1/2)$. The translation curve segments $t_{A_1A_2}$, $t_{A_2A_3}$, $t_{A_3A_1}$ lie
on the coordinate plane $[y,z]$ and the interior angle sum of this translation triangle is $\sum_{i=1}^3(\omega_i)=\pi$.}
\label{}
\end{figure}
The endpoints $T_i^j$ of the position vectors $\bt_i^j=\overrightarrow{E_0T_i^j}$
lie on the unit sphere centred at the origin. The measure of angle $\omega_i$ $(i\in \{1,2,3\})$ of the vectors $\bt_i^j$ and $\bt_r^s$ is equal to the spherical
distance of the corresponding points $T_i^j$ and $T_r^s$ on the unit sphere (see Fig.~3). Moreover, a direct consequence of equations (3.9) that each point pair
($T_2$, $T_1^2$), $(T_3$,$T_1^3$), ($T_2^3$,$T_3^2$)
contains antipodal points related to the unit sphere with centre $E_0$.

Due to the antipodality $\omega_1=T_2E_0T_3 \angle =T_1^2E_0T_1^3 \angle$, therefore their corresponding spherical
distances are equal, as well (see Fig.~3).
Now, the sum of the interior angles $\sum_{i=1}^3(\omega_i)$ can be considered as three consecutive spherical arcs $(T_3^2 T_1^2)$, $(T_1^2 T_1^3)$,
$T_1^3 T_2^3)$.
Since the triangle inequality holds on the sphere, the sum of these arc lengths is greater or equal to the half
of the circumference of the main circle on the unit sphere {i.e.} $\pi$. $\square$
\medbreak
The following lemma is an immediate consequence of the above proof:
\begin{lem}
The angle sum $\sum_{i=1}^3(\omega_i)$ of a $\SOL$ translation triangle $A_1A_2A_3$ is $\pi$ if and only if the points $T_i^j$
$((i,j) \in \{(1,3),(1,2),$ $(2,3),(3,2),(3,0),(2,0)\})$ lie in an Euclidean plane (Fig.~4).
\end{lem}
\begin{lem}
If the vertices of a translation triangle $A_1A_2A_3$ lie in a cooordinate plane of the model of $\SOL$ geometry (see Section 2)
or in a plane parallel to a coordinate plane then the interior angle sum $\sum_{i=1}^3(\omega_i)=\pi$.
\end{lem}
{\bf{Proof}}: We get from equation (2.12) of the translation curves that a point $P$ lies in a coordinate plane
then the corresponding tranlation curve $t_{E_0P}$ also lies in the same coordinate plane.

Moreover, a direct consequence of formulas (2.3) and (2.6) than if a translation triangle $A_1A_2A_3$ lies in a coordinate plane $\alpha$
then its translated image by an orthogonal translation to $\alpha$
is in a to $\alpha$ parallel plane and each to $\alpha$ parallel plane can be derived as a tranlated copy of $\alpha$. $\square$
\medbreak
We can determine the interior angle sum of arbitrary translation triangle.
In the following table we summarize some numerical data of interior angles of given transaltion triangles:
\medbreak
\centerline{\vbox{
\halign{\strut\vrule\quad \hfil $#$ \hfil\quad\vrule
&\quad \hfil $#$ \hfil\quad\vrule &\quad \hfil $#$ \hfil\quad\vrule &\quad \hfil $#$ \hfil\quad\vrule &\quad \hfil $#$ \hfil\quad\vrule
\cr
\noalign{\hrule}
\multispan5{\strut\vrule\hfill\bf Table 1: ~ $A_2(1,-1,1,1),$ ~ $A_3(1,1/2,5,z^3)$ \hfill\vrule}%
\cr
\noalign{\hrule}
\noalign{\vskip2pt}
\noalign{\hrule}
z^3 & \omega_1 & \omega_2 & \omega_3  & \sum_{i=1}^3(\omega_i)  \cr
\noalign{\hrule}
-10 & 1.378505 & 1.52957 & 0.39949 & 3.30757 \cr
\noalign{\hrule}
-2 & 1.37467 & 1.45044 & 0.41389 &  3.23900 \cr
\noalign{\hrule}
-1 & 1.36841 & 1.31743 & 0.48434 &  3.17018 \cr
\noalign{\hrule}
1/100 & 1.35376 & 1.04468 & 0.74818 & 3.14661 \cr
\noalign{\hrule}
1/10 & 1.35196 & 1.01850 & 0.77962  & 3.15008 \cr
\noalign{\hrule}
1/2 & 1.34369 & 0.91985 & 0.90711  & 3.17066 \cr
\noalign{\hrule}
3/4 & 1.33931 & 0.87828 & 0.96332  & 3.18092 \cr
\noalign{\hrule}
3/2 & 1.34516 & 0.83131 & 0.98842  & 3.16489 \cr
\noalign{\hrule}
2 & 1.37178  & 0.83021 &  0.94235  & 3.14433 \cr
\noalign{\hrule}
5 & 1.46886 & 0.84547 & 0.86833 &  3.18265 \cr
\noalign{\hrule}
10 & 1.47522 & 0.84678 & 0.86665 &  3.18866 \cr
\noalign{\hrule}
}}}
\medbreak
\medbreak
\centerline{\vbox{
\halign{\strut\vrule\quad \hfil $#$ \hfil\quad\vrule
&\quad \hfil $#$ \hfil\quad\vrule &\quad \hfil $#$ \hfil\quad\vrule &\quad \hfil $#$ \hfil\quad\vrule &\quad \hfil $#$ \hfil\quad\vrule
\cr
\noalign{\hrule}
\multispan5{\strut\vrule\hfill\bf Table 2: ~ $A_2(1,-1,1,1),$ ~ $A_3(1,1/2,y^3,1/2)$ \hfill\vrule}%
\cr
\noalign{\hrule}
\noalign{\vskip2pt}
\noalign{\hrule}
y^3 & \omega_1 & \omega_2 & \omega_3  & \sum_{i=1}^3(\omega_i)  \cr
\noalign{\hrule}
-10 & 1.90559 & 0.77539 & 0.48862 & 3.16960 \cr
\noalign{\hrule}
-2 & 1.99438 & 0.39617 & 0.86884 &  3.25939 \cr
\noalign{\hrule}
-1 & 2.02152 & 0.38864 & 0.84198 &  3.25214 \cr
\noalign{\hrule}
1/100 & 1.89224 & 0.42533 & 0.83598 & 3.15355 \cr
\noalign{\hrule}
1/10 & 1.86415 & 0.43075 & 0.85319  & 3.14808 \cr
\noalign{\hrule}
1/2 & 1.73149 & 0.45855 & 0.95244  & 3.14248 \cr
\noalign{\hrule}
3/4 & 1.65752 & 0.47867 & 1.01153  & 3.14772 \cr
\noalign{\hrule}
3/2 & 1.51011 & 0.54873 & 1.10619  & 3.16502 \cr
\noalign{\hrule}
2 & 1.45565  & 0.60090 &  1.11440  & 3.17095 \cr
\noalign{\hrule}
5 & 1.34369 & 0.91985 & 0.90711  & 3.17066 \cr
\noalign{\hrule}
10 & 1.30564 & 1.27407 & 0.58095 &  3.16067 \cr
\noalign{\hrule}
}}}
\medbreak
%

%{\bf{Acknowledgement:}}

\end{document}